\input amstex
\documentstyle{amsppt}
%
\catcode`@=11
\redefine\output@{%
  \def\break{\penalty-\@M}\let\par\endgraf
  \ifodd\pageno\global\hoffset=105pt\else\global\hoffset=8pt\fi  
  \shipout\vbox{%
    \ifplain@
      \let\makeheadline\relax \let\makefootline\relax
    \else
      \iffirstpage@ \global\firstpage@false
        \let\rightheadline\frheadline
        \let\leftheadline\flheadline
      \else
        \ifrunheads@ 
        \else \let\makeheadline\relax
        \fi
      \fi
    \fi
    \makeheadline \pagebody \makefootline}%
  \advancepageno \ifnum\outputpenalty>-\@MM\else\dosupereject\fi
}
\def\Beta{\mathchar"0\hexnumber@\rmfam 42}
\catcode`\@=\active
\nopagenumbers
\def\negskp{\hskip -2pt}


\def\blue#1{#1}

\catcode`#=11\def\diez{#}\catcode`#=6
\catcode`_=11\catcode`_=8
\def\mycite#1{\cite{\blue{#1}}\immediate\special{ps:
     ShrHPSdict begin /ShrBORDERthickness 0 def}}

\def\mytag#1{%
    \tag#1}
\def\mythetag#1{\thetag{\blue{#1}}\immediate\special{ps:
     ShrHPSdict begin /ShrBORDERthickness 0 def}}
\def\myrefno#1{\no#1}
\def\myhref#1#2{\blue{#2}\immediate\special{ps:
     ShrHPSdict begin /ShrBORDERthickness 0 def}}

\def\mytheorem#1{\csname proclaim\endcsname{Theorem #1}}
\def\mytheoremwithtitle#1#2{\csname proclaim\endcsname{Theorem #1#2}}

\def\mylemma#1{\csname proclaim\endcsname{Lemma #1}}
\def\mylemmawithtitle#1#2{\csname proclaim\endcsname{Lemma #1#2}}

\def\mycorollary#1{\csname proclaim\endcsname{Corollary #1}}

\def\myconjecture#1{\csname proclaim\endcsname{Conjecture #1}}
\def\myconjecturewithtitle#1#2{\csname proclaim\endcsname{Conjecture #1#2}}
\def\mytheconjecture#1{\blue{#1}\immediate\special{ps:
     ShrHPSdict begin /ShrBORDERthickness 0 def}}

\font\eightcyr=wncyr8
\font\tencyr=wncyr10
\pagewidth{360pt}
\pageheight{606pt}
\topmatter
\title
A conjecture on some estimates for integrals
\endtitle
\author
B.\,N.~Khabibullin
\endauthor
\address Bashkir State University, 32, Z.\,Validi Str., Ufa, 450074, Russia\newline
\vphantom{a}\kern 12pt Cell Phone: +7(347)2736718
\endaddress
\email \myhref{mailto:Khabib-Bulat\@mail.ru}{Khabib-Bulat\@mail.ru}
\endemail
\urladdr
\vtop to 20pt{\hsize=280pt\noindent
\myhref{http://math.bsunet.ru/khb}
{http:/\negskp/math.bsunet.ru/khb}\vss}
\endurladdr
\abstract
    A conjecture concerning some pairs of interfering estimates for some integrals is
formulated in three equivalent versions. Its importance for the the Paley problem for
plurisubharmonic functions and for certain classes of extremal problems for entire 
functions of several variables is declared. 
\endabstract
\endtopmatter
\TagsOnRight
\document

     A {\tencyr\char '074}positive function{\tencyr\char '076} below means a function
whose values are greater than or equal to zero. A real valued function $\phi$ of a real
argument is called an {\tencyr\char '074}increasing function{\tencyr\char '076} if
$x_1\leqslant x_2$ implies the non-strict inequality $\phi(x_1)\leqslant\phi(x_2)$.\par
     The following conjecture associated with solving the Paley problem for
plurisubharmonic functions and with some extremal problems in the theory of entire 
functions of several variables was suggested in \mycite{1}. It is an open problem since 
approximately 1992.

\myconjecture{1} 
Let $S$ be a positive increasing function 
on the ray $[0,+\infty)$ such that $S(0)=0$ and, moreover, let $S$ be logarithmically convex, 
i\.\,e\. the function $x\mapsto S(e^x)$ is convex on the interval $[-\infty,+\infty)$. Let
$\lambda\geqslant 1/2$, $n\geqslant 2$, $n\in\Bbb N$. Under these assumptions if
$$
\gather
\hskip -2em
\int\limits^{\,1}_0 S(tx)\,(1-x^2)^{n-2}\,x\,dx\leqslant t^\lambda\text{\ \ for all \ }
0\leqslant t<+\infty,
\mytag{1}\\
\vspace{-2ex}
\intertext{then}
\vspace{-2ex}
\hskip -2em
\int\limits^{+\infty}_0 S(t)\frac{t^{2\,\lambda-1}}{(1+t^{2\,\lambda})^2}\,dt
\leqslant\frac{\pi\,(n-1)}{2\,\lambda}\prod^{n-1}_{k=1}\Bigl(1+\frac{\lambda}{2\,k}
\Bigr).
\mytag{2}
\endgather
$$
\endproclaim
      Note that the product in the right hand side of \mythetag{2} is associated
with Euler's Beta and Gamma functions usually demoted through $\Beta$ and $\Gamma$:
$$
\hskip -2em
\gathered
\Beta(\lambda/2,n)=\int\limits^{\,1}_0 x^{\lambda/2-1}\,(1-x)^{n-1}\,dx=
\frac{\Gamma(\lambda/2)\ \Gamma(n)}{\Gamma(\lambda/2+n)}=\\
=\frac{\Gamma(\lambda/2)\ (n-1)!}{\Gamma(\lambda/2)\ (\lambda/2)\ (\lambda/2+1)\cdot\ldots\cdot(\lambda/2+(n-1))}=\frac{2}{\lambda}\,\frac{1}
{\dsize\prod^{n-1}_{k=1}\Bigl(1+\frac{\lambda}{2\,k}
\Bigr)},
\endgathered
\mytag{3}
$$
i\.\,e\. the right hand side of \mythetag{2} can be written as 
$$
\frac{\pi (n-1)}{2\,\lambda}\,\prod^{n-1}_{k=1}\Bigl(1+\frac{\lambda}{2\,k}
\Bigr)=\frac{\pi (n-1)}{2\,\lambda^2}\cdot\frac{1}{B(\lambda/2,n)}.
$$
If we choose 
$$
S(t)=2\,(n-1)\prod^{n-1}_{k=1}\Bigl(1+\frac{\lambda}{2\,k}\Bigr)t^\lambda,\qquad t\geqslant 0,\qquad\lambda\geqslant\frac{1}{2}, 
$$
then the inequalities \mythetag{1} and \mythetag{2} turn to the equalities.\par
     The conjecture~\mytheconjecture{1} is valid for $\lambda\leqslant 1$ even without
conditions like convexity for $S$. It is possible to show, though it is difficult, that
for $\lambda>1$ one inevitably should impose some conditions of this sort for $S$.\par
     Beyond these facts, for example, it is even unknown if the conjecture~\mytheconjecture{1}
is valid in the case where $n=2$ for at least one value of $\lambda>1$.\par
     Below we consider the case $\lambda>1$ only.\par
     In explicit form the conjecture~\mytheconjecture{1} is formulated in the end of our
review \mycite{2}.
\mylemmawithtitle{1}{ \rm (was formulated as the Proposition 5.1 in \mycite{3})} A real
valued function $S=S(x)$ on the ray $[0,+\infty)$ such that $S(0)=0$ is an increasing
logarithmically convex function if and only if there is an increasing function $s=s(t)$
on the ray $[0,+\infty)$ such that $S$ is presented as
$$
S(x)=\int\limits^{\kern 1.5pt x}_0\frac{s(t)}{t}\,dt.
$$
\endproclaim
      Using this lemma, we can write the integral from \mythetag{1} as
$$
\int\limits^{\,1}_0 S(t)\,\frac{t^{2\,\lambda-1}}{(1+t^{2\,\lambda})^2}\,dt
=-\frac{1}{2\,(n-1)}\int\limits^{\,1}_0\!\left(\,\int^{tx}_0\frac{s(\tau)}{\tau}\,d\tau
\right)d(1-x^2)^{n-1}.
$$
Integrating this equality by parts, we get
$$
\int\limits^{\,1}_0 S(t)\,\frac{t^{2\,\lambda-1}}{(1+t^{2\,\lambda})^2}\,dt
=-\frac{1}{2\,(n-1)}\int\limits^{\,1}_0\frac{s(tx)}{x}\,(1-x^2)^{n-1}\,dx.
$$
Similarly for the integral \mythetag{2}, we have 
$$
\int\limits^{+\infty}_0\!\!S(t)\frac{t^{2\lambda-1}}{(1+t^{2\lambda})^2}\,dt 
=\frac{1}{2\lambda}\int\limits^{+\infty}_0\frac{s(t)}{t}\,\frac{dt}{1+t^{2\,\lambda}}\,.
$$
Thus the inequalities \mythetag{1} and \mythetag{2} are transformed to the following
inequalities
$$
\gather
\frac{1}{2\,(n-1)}\int\limits^{\,1}_0\frac{s(tx)}{x}\,(1-x^2)^{n-1}\,dx\leqslant t^\lambda
\text{\ \ for all \ }0\leqslant t<+\infty,\\
\frac{1}{2\lambda}\int\limits^{+\infty}_0\frac{s(t)}{t}\,\frac{dt}{1+t^{2\,\lambda}}
\leqslant\frac{\pi\,(n-1)}{2\,\lambda}\prod^{n-1}_{k=1}\Bigl(1+\frac{\lambda}{2\,k}
\Bigr)=\frac{\pi (n-1)}{2\,\lambda^2}\cdot\frac{1}{B(\lambda/2,n)},
\endgather
$$
where $s\geqslant 0$ is an increasing function and the identity \mythetag{3} is used
in deriving the above relationships. Note that without loss of generality we can replace 
$s$ by an increasing function $h\geqslant 0$ defined though the formula
$$
h(x^2)=\frac{1}{4(n-1)}\,s(x)\text{\ \ for \ }x\geqslant 0.
$$
As a result the above two relationships are transformed to a condition for the function 
$h$ and to an inequality for this function:
$$
\gather
2\int\limits^{\,1}_0\frac{h(t^2\,x^2)}{x}\,(1-x^2)^{n-1}\,dx\leqslant (t^2)^{\lambda/2}
\text{\ \ for all \ }0\leqslant t<+\infty,\\
2\!\int\limits^{+\infty}_0\frac{h(t^2)}{t}\,\frac{dt}{1+t^{2\,\lambda}}
\leqslant\frac{\pi}{2}\prod^{n-1}_{k=1}\Bigl(1+\frac{\lambda}{2\,k}
\Bigr)=\frac{\pi}{\lambda}\cdot\frac{1}{B(\lambda/2,n)}.
\endgather
$$
Now, upon the following changes of variables 
$$
\xalignat 3
&x^2=x',
&&t^2=t'
&&\lambda/2=\alpha>1/2
\endxalignat
$$
and redesignating the variables $x'$ and $t'$ back through $x$ and $t$ we get
$$
\gather
\hskip -2em
\int\limits^{\,1}_0\frac{h(t\,x)}{x}\,(1-x)^{n-1}\,dx\leqslant t^\alpha
\text{\ \ for all \ }0\leqslant t<+\infty,
\mytag{6a}\\
\int\limits^{+\infty}_0\frac{h(t)}{t}\,\frac{dt}{1+t^{2\,\alpha}}
\leqslant\frac{\pi}{2}\prod^{n-1}_{k=1}\Bigl(1+\frac{\alpha}{k}
\Bigr)=\frac{\pi}{2\,\alpha}\cdot\frac{1}{B(\alpha,n)}.
\mytag{6b}
\endgather
$$
Thus the conjecture~\mytheconjecture{1} with $\lambda>1$ is equivalent to the
following more simple conjecture.
\myconjecture{2} Let $\alpha>1/2$. Then for any increasing function $h\geqslant 0$
on the ray $[0,+\infty)$ the condition \mythetag{6a} implies the inequality
\mythetag{6b}.
\endproclaim
     There is one more version of the conjectures~\mytheconjecture{1} and 
\mytheconjecture{2}. First of all note that for some rather evident reasons it 
is sufficient to prove these conjectures for smooth functions. Therefore we
can assume that the function $h$ in the conjecture~\mytheconjecture{2} has a
continuous derivative, i\.\,e\. $q=h'\geqslant 0$ on the interval $(0,+\infty)$.
Then, integrating by parts, we get a conjecture equivalent to the previous
conjectures~\mytheconjecture{1} and \mytheconjecture{2}.
\myconjecture{3} Let $\alpha>1/2$. If $q$ is a positive continuous function
on the ray $[0,+\infty)$, then the condition 
$$
\hskip -2em
\int\limits^{\,1}_0\left(\,\,\int\limits^{\,1}_x(1-y)^{n-1}\,\frac{dy}{y}\right)
q(tx)\,dx\leqslant t^{\alpha-1}
\mytag{7}
$$
implies the estimate 
$$
\hskip -2em
\int\limits^{+\infty}_0 q(t)\,\ln\Bigl(1+\frac{1}{t^{\,2\kern 0.2pt\alpha}}\Bigr)\,dt
\leqslant\pi\,\alpha\prod^{n-1}_{k=1}\Bigl(1+\frac{\alpha}{k}\Bigr).
\mytag{8}
$$
\endproclaim
\head
Acknowledgments.
\endhead
I am grateful to Dr. Ruslan Sharipov for translating this paper into English and for
typesetting the English version of this paper.


\Refs
\ref\myrefno{1}\by Khabibullin~B.~N.\paper Paley problem for plurisubharmonic functions 
of a finite lower order\jour Sbornik: Mathematics\vol 190\issue 2\yr 1999\pages 309--321 
\endref
\ref\myrefno{2}\by Khabibullin~B.~N.\paper The representation of a meromorphic function as 
a quotient of entire functions and the Paley problem in $\Bbb C^n$: survey of some results
\jour Mathematical Physics, Analysis, and geometry (Ukraine) \yr 2002\vol 9\issue 2
\pages 146--167\moreref see also
\myhref{http://arxiv.org/abs/math.CV/0502433}{Archiv:math.CV/0502433}
\endref
\ref\myrefno{3}\by Kondratyuk~A.~A.\book Fourier series and meromorphic functions
\publ {\eightcyr\char '074 Vishcha shkola\char '076} publishers
\publaddr Lviv (Ukraine)\yr 1988
\endref
\endRefs
\enddocument
\end